\documentclass[11pt]{article}
\usepackage{url}
\begin{document}
\baselineskip=15pt
\title{A Comment on Matiyasevich's Identity \#0102 with Bernoulli Numbers}
\author{H. Gopalkrishna Gadiyar and R. Padma\\ ~~\\AU-KBC Research Centre\\ M. I. T. Campus of Anna University\\ Chromepet, Chennai 600 044, INDIA\\E-mail: \{gadiyar, padma\}@au-kbc.org}
\date{~~~}
\maketitle

\begin{abstract}We connect and generalize Matiyasevich's identity \#0102 with Bernoulli numbers and an identity of Candelpergher, Coppo and Delabaere on Ramanujan summation of the divergent series of the infinite sum of the harmonic numbers. The formulae are analytic continuation of Euler sums and lead to new recursion relations for derivatives of Bernoulli numbers. The techniques used are contour integration, generating functions and divergent series.
\end{abstract}

\section{Introduction} In his personal journal on the Internet Matiyasevich \cite{matiyasevich} has presented the identity \#0102 with Bernoulli numbers: \emph{for even $n~>2$},
\begin{equation}
\sum_{k+l=n} \frac{B_k}{k}B_l ~-~\sum_{k+l=n} \left(\!\!\!
  \begin{array}{c}
    n \\
    k
  \end{array}
  \!\!\!\right) \frac{B_k}{k}B_l -B_n H_n =0 \, , \label{eq:mat}
\end{equation}
and welcomed readers of his website to volunteer further information on the topic. Here, $B_n$'s are the Bernoulli numbers given by
\begin{equation}
\frac{z}{e^z-1}=\sum_{n=0}^\infty (-1)^n~B_n ~\frac{z^n}{n!} \label{eq:ber}
\end{equation}
and $H_n=\sum_{j=1}^{n}\frac{1}{j}$ is the $n^{th}$ harmonic number. Note that the definition of Bernoulli numbers given above \cite{woon} changes only the value of $B_1$ from the usual definition and Matiyasevich's identity is not affected by this modification. In this paper, we use the well known tools of generating functions \cite{wilf}, contour integration and divergent series \cite{hardy} to connect Matiyasevich's discovery to two other classes of identities. The first one due to Candelpergher, Coppo and Delabaere \cite{Candelpergher:Coppo:Delabaere} is given by
\begin{equation}
\sum_{n \ge 1}^{[\cal R]}H_n = \frac{3}{2}\gamma +\frac{1}{2}-\frac{1}{2}\log (2\pi) \, , \label{eq:candel}
\end{equation}
where $[\cal R]$ represents the Ramanujan summation and $\gamma $ is  Euler-Mascheroni constant. 

The other class of identities related to Matiyasevich's discovery is given by explicit formulae for Euler sums in terms of Riemann-zeta values. We give a sample of typical examples of such identities below.
\begin{eqnarray}
\sum_{n=1}^\infty \frac{H_n}{n^2}&=&2 \zeta (3)\\
\sum_{n=1}^\infty \frac{H_n}{n^3}&=&\frac{5}{4}\zeta (4)\\
\sum_{n=1}^\infty \frac{H_n}{n^4}&=&3 \zeta (5) - \zeta (2) \zeta (3)
\end{eqnarray}
The left hand side of these identities are examples of Euler sums. For the definition of general class of Euler sums and their formulae, see \cite{Flajolet:Salvy}. Apostol and Vu \cite{apostol:vu} and Matsuoka \cite{matsuoka} studied the related Dirichlet series \( h(s)=\sum_{n=1}^\infty \frac{H_n}{n^s} \), $\Re (s) >1$ and showed that $h(s)$ can be analytically continued to the entire complex plane. For the history of Euler sums, see \cite{berndt} and for a recent survey see \cite{bowman:bradley} and the references therein. Proof and generalizations of Matiyasevich's identity and its relation to Miki's identity are discussed in \cite{pan:sun}, \cite{shiratani:yokoyama}, \cite{gessel} and \cite{dunne:schubert}. 

\section{A Hybrid of Matiyasevich and Candelpergher-Coppo-Delabaere Identities}
\label{sect:hybrid}

In this section we prove the following result which combines the results of Matiyasevich and Candelpergher-Coppo-Delabaere.

\noindent For $n \ge 1$,
\begin{equation} 
(-1)^{n-1}~n~h(-n+1)=B'_n+n B_{n-1}+\gamma B_n- \sum_{k+l=n} \left(\!\!\!
  \begin{array}{c}
    n \\
    k
  \end{array}
  \!\!\!\right) \frac{B_k}{k}B_l -B_n H_n \label{eq:matext}
\end{equation}
where \( h(-k)= \sum_{n=1}^\infty H_n n^k \) is understood in the sense of an appropriate sum of a divergent series.

\noindent {\bf Remarks} 1. The last two terms are identical to those in Matiyasevich's identity (\ref{eq:mat}).\\ 2. The first identity in this set is that of Canderpergher,Coppo and Delabaere. By substituting $n=1$ and using the fact \cite{woon} that $\displaystyle{B_1'=\frac{1}{2}-\frac{1}{2}\log 2\pi}$ we get
\begin{equation} 
h(0)=\sum_{n = 1}^{\infty}H_n = \frac{1}{2}\gamma +\frac{1}{2}-\frac{1}{2}\log (2\pi) \, . \label{eq:ours}
\end{equation}
The difference between the identities (\ref{eq:candel}) and (\ref{eq:ours}) comes from the fact that the Ramanujan sum of the divergent series $\displaystyle{\sum_{n\ge 1}^{[\cal R]} 1~=~\frac{1}{2}}$ whereas the value arrived at by standard summation methods is $\displaystyle{-\frac{1}{2}}$. Such differences often happen when differing methods of summing divergent series are used. See for example Chapter 5 of \cite{soule:abramovich:burnol:kramer} for a recent exposition.
 
\noindent {\bf Proof.} Let us take the term 
\begin{equation}
\sum_{k+l=n} \left(\!\!\!
  \begin{array}{c}
    n \\
    k
  \end{array}
  \!\!\!\right) \frac{B_k}{k}B_l =n!~\sum_{k+l=n}\frac{B_k}{k~k!}\frac{B_l}{l!} \label{eq:mat1}
\end{equation}
in Matiyasevich's identity (\ref{eq:mat}). Recall that \cite{wilf}, if $\displaystyle{a(z)=\sum \frac{a(n)}{n!} z^n}$ and $\displaystyle{b(z)=\sum \frac{b(n)}{n!} z^n}$, then $\displaystyle{a(z) b(z)=\sum \frac{c(n)}{n!} z^n}$, where $\displaystyle{\frac{c(n)}{n!}=\sum_{k+l=n} \frac{a(k)}{k!}\frac{b(l)}{l!}}$. Choosing the generating functions 
\begin{equation}
\frac{-z}{e^{-z}-1}=\sum_{n=0}^\infty B_n ~\frac{z^n}{n!} \, , \label{eq:ber1}
\end{equation}
(set $z \rightarrow -z$ in (\ref{eq:ber})) and
\begin{equation}
\log (e^z-1)-\log z=\log \left (\frac{e^z-1}{z}\right ) = \sum_{n=1}^\infty \frac{B_n}{n}~\frac{z^n}{n!} \, ,
\end{equation}
see \cite{gessel} and \cite{berndt}, one gets (\ref{eq:mat1}). That is,
\begin{eqnarray}
\sum_{n=1}^\infty \sum_{k+l=n}\frac{B_k}{k~k!}\frac{B_l}{l!} z^n&=&\left [\log \left (\frac{e^z-1}{z}\right ) \right ] \left [\frac{-z}{e^{-z}-1} \right ]\nonumber \\
&=&\left (\log \left (e^z-1 \right ) - \log z \right ) \left (\frac{-z}{e^{-z}-1} \right ) \nonumber \\
~&=& (-z) \frac{\log(e^z-1)}{e^{-z}-1}  ~+ ~\frac{z}{e^{-z}-1}~ \log z\nonumber \\
~&=& (-z) \frac{\log e^z(1-e^{-z})}{e^{-z}-1}  ~+~\frac{z}{e^{-z}-1}~\log z\nonumber \\
~&=& \frac{-z^2}{e^{-z}-1} +z \frac{\log (1-e^{-z})}{1-e^{-z}}~+~ \frac{ze^z}{1-e^{z}} ~ \log z\, . \nonumber\\
~& & ~~~~~~~~~~~~~~~~~~~~~\label{eq:newid}
\end{eqnarray}
In order to pick out the $n^{th}$ coefficient of the series in (\ref{eq:newid}), let us integrate along the Hankel contour $C$. $C$ is a loop that starts at infinity on the negative real axis, encircles the origin in the positive direction, excluding the points $\pm 2i\pi$, $\pm 4i\pi \cdots$, and returns to negative infinity. Thus,
\begin{equation}
\sum_{k+l=n}\frac{B_k}{k~k!}\frac{B_l}{l!} =\frac{1}{2\pi i} \int_C \frac{1}{z^{n+1}} \sum_{n=1}^\infty \sum_{k+l=n}\frac{B_k}{k~k!}\frac{B_l}{l!} z^n dz \, ,
\end{equation}
and from (\ref{eq:newid}), we get
\begin{eqnarray}
\sum_{k+l=n}\frac{B_k}{k~k!}\frac{B_l}{l!} &=&\frac{1}{2\pi i} \int_C \frac{1}{z^{n}} \frac{-z}{e^{-z}-1}~ dz \nonumber\\
~& &~~ + \frac{1}{2\pi i} \int_C \frac{1}{z^{n}}  \frac{\log (1-e^{-z})}{1-e^{-z}}~dz \nonumber\\
~& &~~ + \frac{1}{2\pi i} \int_C \frac{1}{z^{n}}~ \frac{e^z}{1-e^{z}}~\log z ~dz \label{eq:int}\\
~&=& \frac{B_{n-1}}{(n-1)!} - \frac{(-1)^{n-1}h(-n+1)}{(n-1)!} +\frac{B'_n}{n!} \nonumber\\
~&&~~-\frac{B_n}{n!} \left (-\gamma +H_n\right ) \label{eq:bernoulli} \, ,
\end{eqnarray}
where \( B'_n = \frac{d}{ds}B_s|_{s=n}\) is the derivative of the Bernoulli number. Multiplying both sides of (\ref{eq:bernoulli}) by $n!$ gives (\ref{eq:matext}) which is an extension of Matiyasevich's identity \#0102. The evaluation of the three integrals is given below as Lemmas 1, 2 and 3.
\newtheorem{lemma}{Lemma} 
\begin{lemma}
For $n \ge 1$,
\begin{equation}
\frac{1}{2\pi i} \int_C \frac{1}{z^{n}} \frac{-z}{e^{-z}-1}~ dz = \frac{B_{n-1}}{(n-1)!} \, .
\end{equation} 
\end{lemma}
\noindent{\bf Proof.} Immediate from (\ref{eq:ber1}).

\begin{lemma} For $n \ge 1$
\begin{equation}
\frac{1}{2\pi i} \int_C \frac{1}{z^{n}}  \frac{\log (1-e^{-z})}{1-e^{-z}}~dz = \frac{(-1)^n h(-n+1)}{(n-1)!} \, . 
\end{equation}
\end{lemma}
\noindent {\bf Proof}  The integrals along the upper and lower lips of the Hankel contour $C$ cancel each other and thus the integral is equal to the residue of $\displaystyle{\frac{1}{z^{n}}  \frac{\log (1-e^{-z})}{1-e^{-z}}}$ at the origin. We now evaluate this residue. Using the identity
\begin{equation}
-\frac{\log(1-x)}{1-x}~=~ \sum_{n=1}^\infty H_n x^n \, ,
\end{equation}
where $H_n$ is the $n^{th}$ harmonic number defined earlier, we get
\begin{equation}
-\frac{\log (1-e^{-z})}{1-e^{-z}} = \sum_{n=1}^\infty H_n e^{-nz} \, .
\end{equation}
Expanding the exponentials as Taylor series around the origin, and rearranging the divergent series we get 
\begin{equation}
\sum_{n=1}^\infty H_n e^{-nz} = \sum_{n=1}^\infty \sum_{k=0}^\infty H_n \frac{(-nz)^k }{k!}=\sum_{k=0}^\infty \frac{(-z)^k}{k!} h(-k) \, .
\end{equation}

\begin{lemma}
\noindent For $n \ge 1$
\begin{equation}
\frac{1}{2\pi i} \int_C \frac{1}{z^n}~ \frac{e^z}{1-e^z}~\log z ~dz =
\frac{B'_n}{n!} -\frac{B_n}{n!} \left (-\gamma +H_n\right ) \label{eq:berder} \, .
\end{equation}
\end{lemma}

\noindent {\bf Proof} See Lemma 1 of \cite{boyadzhiev}. We have \cite{apostol}
\begin{equation}
-\frac{B_{s+1}}{\Gamma(s+2)} = \frac{1}{2\pi i} \int_C \frac{1}{z^{s+1}}~ \frac{e^z}{1-e^{z}}~dz \, . \label{eq:interpolation}
\end{equation}
By taking derivatives on both sides of (\ref{eq:interpolation}) with respect to $s$, we get
\begin{equation}
\frac{B'_{s+1}}{\Gamma(s+2)} - \frac{B_{s+1}}{\Gamma(s+2)} \frac{\Gamma '(s+2)}{\Gamma(s+2)}  =\frac{1}{2\pi i} \int_C \frac{1}{z^{s+1}}~ \frac{e^z}{1-e^{z}}~\log z ~dz \, .\label{eq:berderivative}
\end{equation}
Using the identity \cite{Copson}
\begin{equation}
\frac{\Gamma '(n+1)}{\Gamma(n+1)}= -\gamma +H_n \, ,
\end{equation}
and putting $s=n-1$ in (\ref{eq:berderivative}), we get (\ref{eq:berder}).

\section{Possible Directions for Generalization}
\label{sect:generalize}

Let us take the integral 
\begin{eqnarray}
\int_0^\infty \frac{x^{s-1}}{(e^x-1)^2} ~dx &=& \int_0^\infty x^{s-1} \left ( e^{-2x}+2e^{-3x}+3e^{-4x}+\cdots \right )x^{s-1}~dx\nonumber \\
&=& \left ( \frac{1}{2^s}+ \frac{2}{3^s}+\frac{3}{4^s}++\cdots \right )\Gamma (s)\nonumber \\
&=& \left ( \frac{2-1}{2^s}+ \frac{3-1}{3^s}+\frac{4-1}{4^s}++\cdots \right )\Gamma (s)\nonumber \\
&=& \Gamma (s)\left ( \zeta (s-1) - \zeta (s) \right )\, ,
\end{eqnarray}
if $\Re (s) >2$ \cite{titchmarsh}. Writing this real integral as complex integral using the Hankel contour \cite{apostol} gives 
\begin{equation}
\zeta(s-1)-\zeta(s)= \Gamma(1-s) \frac{1}{2\pi i} \int_C \frac{z^{s-1} e^{2z}}{(1-e^z)^2} ~dz \label{eq:Js}
\end{equation}
The integral on the right hand side can be analytically continued to the entire complex plane as is done in \cite{apostol}, and now replacing $s$ by $-s$ in (\ref {eq:Js}) we get
\begin{equation}
\frac{\zeta(-s-1)-\zeta(-s)}{\Gamma(1+s)} = \frac{1}{2\pi i} \int_C \frac{1}{z^{s+1}} \frac{e^{2z}}{(1-e^z)^2} ~dz \, .
\end{equation}
Differentiating  both sides with respect to $s$ and writing $\log z = \log(e^z-1)-\log(\frac{e^z-1}{z})$, the right hand side of the integral at $s=n$ will give,
\begin{eqnarray}
\frac{1}{2\pi i} \int_C \frac{1}{z^{n+1}} \frac{e^{2z} \log z }{(1-e^z)^2} ~dz &=& \frac{1}{2\pi i} \int_C \frac{1}{z^{n+1}}\left (\frac{e^z \log z }{e^z-1} \right ) \left ( \frac{e^z}{e^z-1} \right )~dz \nonumber \\
&=&\frac{1}{2\pi i} \int_C \frac{1}{z^{n+1}} \left (\log(e^z-1)\frac{e^z}{e^z-1}\right )\left ( \frac{e^z}{e^z-1}\right ) dz\nonumber \\ 
&&~-\frac{1}{2\pi i} \int_C \frac{1}{z^{n+1}} \left (\log \left (\frac{e^z-1}{z} \right )\right ) \left ( \frac{e^z}{e^z-1}\right ) dz\, .\nonumber \\
&&~
\end{eqnarray}
Writing the power series for the two terms in terms of $h(-k)$ and $B_n$ will give relationship between $B'_{n+1}$, $B'_n$ and $h(-k)$ for $k \le n$. As we have already derived a connection between $h(-n+1)$ and $B'_n$ in Section \ref{sect:hybrid}, this would lead to a recursion relation for $B'_n$ and Glaisher-Kinkelin-Bendersky constants \cite{glaisher}, \cite{dowker}. It has not escaped our attention that finding a recursion for $B'_n$ will lead to a formula for $\zeta (2k+1) $ from the functional equation of the $\zeta $-function. 

\section{Comments} 
\begin{enumerate}
\item There is a lot of confusion in the literature regarding Riemann zeta functions, Bernoulli numbers and their representations. There are essentially real and complex representations of Bernoulli numbers and Riemann zeta functions. We are inspired by Rota \cite{gessel1} to use the idea that Bernoulli numbers can be generated by a complex representation around the Hankel contour. The idea of using generating functions for Euler sums was powerfully demonstrated by Flajolet. Gessel's derivation of Miki's identities is on similar lines but cancels the logarithmic part. Roman and Rota in \cite{roman:rota} develop the analogue of logarithmic Taylor series and attempt to sum derivatives of Bernoulli numbers. This paper essentially combines the ideas of these authors. 

\item The Norl\"und polynomials $B_n^{(z)}$ are defined by 
\begin{equation}
\sum_{n=0}^\infty B_n^{(z)}\frac{x^n}{n!} = \left ( \frac{x}{e^x-1} \right )^z \, .
\end{equation}
Note that that Matiyasevich's identity (\ref{eq:mat}) is related to the derivative of the above at $z=1$. 

\item As Euler sums are related to zeta functions in several identities as an after thought it seems obvious that using the functional equations for Euler sums and zeta functions, there have to be several identities which are analytically continued versions of existing identities. 

\item In \cite{dunne:schubert} Dunne and Schubert use the generating function
\begin{equation}
\tilde{\psi} (x)=\psi (x) -\log x +\frac{1}{2x} \sim \sum_{k=1}^\infty \frac{B_{2k}}{2k} \frac{1}{x^{2k}} \, ,
\end{equation}
where $\psi (x) = \frac{\Gamma '(x)}{\Gamma (x)}$. If the logarithm term is not killed, using the ideas in \cite{Flajolet:Salvy}, new identities may be obtained by a new method. These identities should also be related to analytic continuation of Euler sums.

\item In \cite{dowker} the importance of the combination  $\zeta '(-n) -H_n \zeta(-n)$ in the theory of Glaisher-Kinkelin-Berndersky constants is stressed. This combination occurs naturally in this derivation.

\item A novel feature of the approach outlined is the surprising appearance of a divergent series as a residue. These formulae cannot be checked easily numerically as they involve divergent series. Further there seems to be some difference between the conclusions of Apostol and Vu \cite{apostol:vu} and Candelpergher, Coppo and Delabaere \cite{Candelpergher:Coppo:Delabaere} regarding the values of the analytically continued Euler sums at negative integers. This would require very rigorous and careful analysis. Hence the derivations we have given should be considered heuristic and merely sketching the rich possibilities existing in this direction.

\end{enumerate}

\section{Acknowledgements} The authors wish to thank Professor H. S. Sharatchandra and Professor M. S. Rangachari for several fruitful discussions in initial stages of the work.

\end{document}